\documentclass[12pt]{article}

\usepackage[cmtip,arrow]{xy}
\usepackage{amsmath,amssymb,enumerate,pb-diagram,pb-xy}

\usepackage[applemac]{inputenc}

\def \1{{\bf 1}}
\def \A{{\mathbb A}}

\def \al{\alpha}

\def \bs{\backslash}
\def \C{{\mathbb C}}

\def \CF{{\cal F}}
\def \CH{{\cal H}}
\def \CL{{\cal L}}
\def \CM{{\cal M}}
\def \CN{{\cal N}}

\def \diag{\operatorname{diag}}
\def \df{\ \begin{array}{c} _{\rm def}\\ ^{\displaystyle =}\end{array}\ }

\def \e{\emph}
\def \eps{\varepsilon}

\def \fin{{\rm fin}}

\def \ga{\gamma}
\def \Ga{\Gamma}
\def \GL{{\rm GL}}
\def \H{{\mathbb H}}
\def \Hom{{\rm Hom}}

\def \Im{{\rm Im}}
\def \k{{\mathfrak k}}

\def \La{\Lambda}

\def \N{{\mathbb N}}

\def \Per{\operatorname{Per}}

\def \prf{{\bf Proof: }}
\def \PGL{\operatorname{PGL}}

\def \Q{{\mathbb Q}}
\def \qed{\ifmmode\eqno \square
        \else\noproof\vskip 12pt plus 3pt minus 9pt \fi}
\def \noproof{{\unskip\nobreak\hfill\penalty50\hskip2em\hbox{}%
     \nobreak\hfill $\square$\parfillskip=0pt%
     \finalhyphendemerits=0\par}}
\def \R{{\mathbb R}}
\def \Re{\operatorname{Re}}
\def \res{{\rm res}}

\def \SL{{\rm SL}}
\def \SO{{\rm SO}}

\def \ul{\underline}
\def \vol{{\rm vol}}

\def \Z{{\mathbb Z}}
\def \({\left(}
\def \){\right)}
\def \={{\ =\ }}

\newcommand{\tto}[1]{\stackrel{#1}{\longrightarrow}}
\newcommand{\dotcup}{\ensuremath{\mathaccent\cdot\cup}}
\newcommand{\stack}[2]{\genfrac{}{}{0pt}{}{#1} {#2}}
\newcommand{\norm}[1]{\left|\hspace{-1pt}\left| #1\right|\hspace{-1pt}
\right|}

\renewcommand{\sp}[1]{\left\langle #1\right\rangle}
\newcommand{\ol}[1]{\overline{#1}}

\newtheorem{theorem}{Theorem}[subsection]

\newtheorem{lemma}[theorem]{Lemma}
\newtheorem{corollary}[theorem]{Corollary}
\newtheorem{proposition}[theorem]{Proposition}
\newtheorem{exmples}[theorem]{Examples}
\newenvironment{examples}[0]{\begin{exmples}{\ }\\ 	\vspace{-20pt}\nopagebreak[4]
	\begin{itemize}\rm}{\end{itemize}\end{exmples}}
\newtheorem{exmple}[theorem]{Example}
\newenvironment{example}[0]{\begin{exmple}\rm}
	{\end{exmple}}

\begin{document}

\pagestyle{myheadings} \markright{AUTOMORPHIC FORMS OF HIGHER ORDER}

\title{Automorphic forms of higher order}
\author{Anton Deitmar \& Nikolaos Diamantis\\ \ \\
Journal of the London Mathematical Society. 80, 18-34 (2009)}
\date{}
\maketitle

MSC: {\bf 11F12}, 11F25, 11F66, 11F99

{\bf Abstract:}
In this paper a theory of Hecke operators for higher order modular forms 
is established.
The definition of higher order forms is extended beyond 
the realm of parabolic invariants. A canonical inner product is introduced.
The role of representation theoretic methods is clarified and, motivated 
by higher order forms, new convolution products of L-functions are 
introduced.

{\small\tableofcontents}

\section*{Introduction}

Higher-order modular forms have in recent years arisen in various 
contexts and they have been studied as analytic functions \cite{Gold,KZ}. In particular 
their spaces \cite{CDO,DO, DS, DSr} have been investigated. Parallel to 
that, L-functions were attached to higher-order forms \cite{DKMO,F,IM} and some 
of their basic aspects studied.
In \cite{ES}, a cohomology theory for higher invariants is developed, and an Eichler-Shimura isomorphism for higher order forms is established.

The current paper serves the following purposes,
\begin{itemize}
\item to establish a theory of Hecke operators on higher order forms,
\item to extend the definition of higher forms beyond parabolic invariants,
\item to clarify the role of representation theoretic methods in the theory, and
\item to introduce new convolution products of L-functions of higher forms.
\end{itemize}

The first item fills a long-standing gap in the theory of higher order forms by constructing a natural Hecke action.
This is surprising, as there is no adelic counterpart of higher order forms.
The Hecke operators form bounded self-adjoint operators on direct limits of spaces of higher order forms.
It is an on-going long-term project of the authors to better understand the spectral decompositions of Hecke operators.
For the second item, one gets that in the case of higher order forms, the crucial Fourier expansion is replaced by a ``Fourier-Taylor-expansion'' which is introduced in this paper.
For the third item, representation theoretic methods, it turns out that a \emph{intervention of Lie groups}, as Dieudonn\'e terms it, 
is possible in the theory and, in fact, higher-order forms can be incorporated into the same representational context as the classical 
automorphic forms. 
One would thus expect a distinction from classical forms in the \emph{intervention of adeles}. 
Indeed, it turns out 
that there is no intervention of adeles, as there are no higher forms on the adelic level.
The last item in the list, the convolution product, is inspired by the second insofar  
as the L-functions of higher-order forms are special cases of the convolution products.
We show analyticity and functional equation in greater generality.

\section{Higher invariants and Hecke operators}
\subsection{Higher invariants}
Let $R$ be a commutative ring with unit and let 
$\Ga$ be a group.
Let $R[\Ga]$ be the group algebra and $I_\Ga\subset R[\Ga]$ the augmentation ideal.
Note that $I_\Ga$ is a free $R$-module with basis $(\ga-1)_{\ga\in\Ga}$.
For an $R[\Ga]$-module $V$, the usual invariants are the elements of the $R$-module $H^0(\Ga,V)=V^\Ga=\{v\in V: I_\Ga v=0\}$.
We define the set of \e{higher invariants} to be
$$
H_q^0(\Ga,V)\=\{ v\in V: I_\Ga^{q+1} v =0\},
$$
for $q=0,1,2,\dots$, where $I_\Ga^k$ is the $k$-th power of the ideal 
$I_\Ga$. 
For $v\in V$ one has $v\in H^0_{q+1}(\Ga,V)$ if and only if $(\ga -1)v\in H^0_{q}(\Ga,V)$ for every $\ga\in\Ga$.
Note that if $\Ga$ is \e{perfect}, i.e., if $[\Ga,\Ga]=\Ga$, then $I_\Ga^{q+1}=I_\Ga$ for every $q\ge 0$ and so in that case there are no higher invariants except the usual invariants.
This is due to the observation
$$
[a,b]-1=aba^{-1}b^{-1}-1= (ab-ba)(ba)^{-1}=[(a-1)(b-1)-(b-1)(a-1)](ba)^{-1},
$$
which shows that if $\Ga$ is perfect, then $I_\Ga^2=I_\Ga$ and hence inductively, $I_\Ga^{q+1}=I_\Ga$.

\begin{proposition}\label{1.1}
Assume $R=\C$, that $\Ga$ is a compact group, and that $V$ is a Hilbert space on which $\Ga$ acts by a continuous representation.
Then there are no higher invariants except the classical invariants, i.e., one has $H_q^0(\Ga,V)=H^0(\Ga,V)$ for every $q\ge 0$.
\end{proposition}
\prf
Let $q\ge 1$.
Formally set $H_{-1}^0(\Ga,V)=0$.
For every $v\in H^0_q(\Ga,V)$ the map 
$\ga\mapsto (\ga-1)v\in H^0_{q-1}(\Ga,V)/H^0_{q-2}(\Ga,V)$ is a group homomorphism to the 
additive group, as $(\ga\tau-1)\equiv (\ga-1)+(\tau-1)\mod I_\Ga^2$.
Let $\ga\in\Ga$.
Then one has $(\ga^n-1)v=n(\ga-1)v$.
Since $\Ga$ is compact, there is a sequence $n_k\to\infty$ in $\N$ such that $\ga^{n_k}$ converges in $\Ga$.
Hence $(\ga^{n_k}-1)v=n_k(\ga-1)v$ converges, too, and so $(\ga-1)v=0$.
This holds for every $\ga\in\Ga$ and so $(\ga-1)v\in H_{q-2}^0$ which implies $v\in H_{q-1}^0$, so $H_q^0=H_{q-1}^0$ and inductively $H_q^0=H_0^0$.
\qed

\subsection{Hecke operators}
We now come to Hecke operators. For this, let
$(G,\Ga)$ be a \e{Hecke pair}, i.e., $G$ is a group and $\Ga$ is a subgroup such that for every $g\in G$ the set $\Ga g\Ga/\Ga$ is finite.

\begin{examples}
\item A classical example is $G=\GL_2(\Q)$ and $\Ga=\GL_2(\Z)$.
\item Let $G$ be a topological group and let $\Ga$ be a compact open subgroup.
Then $(G,\Ga)$ is a Hecke pair, since the compact set $\Ga g\Ga$ can be covered by finitely many open sets of the form $x\Ga$, $x\in G$.
\end{examples}

The \e{Hecke algebra} $\CH=\CH_R(G,\Ga)$ is the $R$-module of all functions $f:\Ga\bs G/\Ga\to R$ of finite support with the convolution product
$$
f*h(x)\= \sum_{y\in G/\Ga} f(y) h(y^{-1}x).
$$
As an $R$-module, $\CH$ is free with basis $(\1_{\Ga g\Ga})_{g\in G}$, where $\1_A$ is the characteristic function of the set $A$.
Let $V$ be an $R[G]$-module.
The Hecke algebra $\CH$ acts naturally on $V^\Ga=H^0(\Ga,V)$ via
$$
f.v\= \sum_{y\in G/\Ga} f(y)yv.
$$
In particular one has
$$
\1_{\Ga g\Ga}v\= \sum_{j=1}^n g_jv,
$$
where $\Ga g\Ga=\dotcup_{j=1}^ng_j\Ga$.
As $v$ is $\Ga$-invariant, this expression does not depend on the choice of the representatives $g_j$.

For $q\ge 1$ and $v\in{H_q^0}(\Ga,V)$, however, the expression 
$\sum_{j=1}^n g_jv$ \emph{will} in general depend on the choice of the representatives $(g_j)$.
Any other set of representatives is of the form $(g_j\ga_j)$ for some $\ga_j\in\Ga$.
Note that
$$
\sum_j g_j \ga_jv-\sum_jg_jv\= \sum_jg_j \underbrace{(\ga_j-1)v}_{\in {H_{q-1}^0}(\Ga,V).}
$$
The group $\Ga$ permutes the finite set $\Ga g\Ga/\Ga$ by left multiplication.
Let $\Ga(g)\subset\Ga$ be the subgroup of all elements which act trivially on $\Ga g\Ga/\Ga$.
Then $\Ga(g)$ is a finite index normal subgroup of $\Ga$.
Note that
$$
\Ga(g)\= \bigcap_{\ga\in\Ga}\ga(\Ga\cap g\Ga g^{-1})\ga^{-1},
$$
so $\Ga(g)$ is the biggest normal subgroup contained in $\Ga\cap g\Ga g^{-1}$.

\begin{lemma}
\begin{enumerate}[\rm (a)]
\item For $v\in H_q^0(\Ga,V)$ the sum $\sum_{j=1}^ng_jv$ lies in $H_q^0(\Ga(g),V)$.
\item For $v\in H_q^0(\Ga,V)$ and any choice of elements $\ga_j\in\Ga$ the sum \\
$\sum_{j=1}^ng_j(\ga_j-1)v$ lies in $H_{q-1}^0(\Ga(g),V)$.
\end{enumerate}\end{lemma}
\prf
(a) We have to show that the sum is anihilated by
$$
(\sigma_0-1)\cdots (\sigma_{q}-1)
$$
for any given $\sigma_0,\dots,\sigma_q\in\Ga(g)$.
By the definition of $\Ga(g)$ it follows that for every $j$ and every $k$
there exist $\eta_j^k\in\Ga$ such that $\sigma_kg_j=g_j\eta_j^k$.
Hence,
$$
(\sigma_0-1)\cdots (\sigma_{q}-1)\sum_{j=1}^ng_j\=\sum_jg_j(\eta_j^0-1)\cdots (\eta_j^{q}-1)v\= 0.
$$
The proof of (b) is similar.
\qed

The lemma implies that we get a well defined Hecke operator
$$
T_{\Ga g\Ga}: {H_q^0}(\Ga,V)/{H_{q-1}^0}(\Ga,V)\ \to\ 
{H_q^0}(\Ga(g),V)/{H_{q-1}^0}(\Ga(g),V).
$$
We need to extend this construction to finite index subgroups $\Sigma\subset\Ga$.
First note that $(G,\Sigma)$ is a Hecke pair again. We abbreviate
$$
\bar H_q(\Sigma,V)\={H_q^0}(\Sigma,V)/{H_{q-1}^0}(\Sigma,V).
$$
As every finite index subgroup contains a finite index normal subgroup it suffices to assume that $\Sigma$ is normal in $\Ga$.
We define $\Sigma(g)$ to be the kernel of the homomorphism $\Sigma\to\Per(\Ga g\Ga /\Sigma)$.
Note that $\Ga$ is present in the definition, although not in the notation.
Then, $\Sigma(g)=\bigcap_{\ga\in\Ga}\ga(\Sigma\cap g\Sigma g^{-1})\ga^{-1}$, and $\Sigma(g)$ is normal of finite index in $\Ga$.
We define
$$
T_{\Ga g\Ga} : \bar H_q(\Sigma,V)\ \to\ \bar H_q(\Sigma(g),V)
$$
as follows.
Write $\Ga g\Ga$ as a disjoint union of $\Sigma$-cosets $\dotcup_j h_j\Sigma$ and set
$$
T_{\Ga g\Ga} v\=\frac 1{[\Ga:\Sigma]}\sum_j h_j v.
$$
The same reasoning as before shows the well-definedness of $T_{\Ga g\Ga}$.
The factor $\frac 1{[\Ga:\Sigma]}$ will make the Hecke operator compatible with change of groups as follows.
Assume $\Sigma'\subset\Sigma$ is another finite index normal subgroup.
Then $\Sigma'(g)\subset \Sigma(g)$ as well and the inclusions 
$$
\begin{array}{ccc}
H_q^0(\Sigma,V) & \subset & H_q^0(\Sigma',V)\\
\cup&&\cup\\
H_{q-1}^0(\Sigma,V) & \subset & H_{q-1}^0(\Sigma',V)
\end{array}
$$
show that there is a natural restriction homomorphism $\res^\Sigma_{\Sigma'}:\bar H_q(\Sigma, V)\to \bar H_q(\Sigma', V)$.
As the intersection of two finite index subgroups is a finite index subgroup, these spaces form a direct system indexed by the set of all finite index normal subgroups $\Sigma$ of $\Ga$.

\begin{lemma}
For any two finite index normal subgroups $\Sigma'\subset\Sigma$ of $\Ga$ and $g\in G$ the diagram
$$
\begin{diagram}
\node{\bar H_q(\Sigma, V)}\arrow{e,t}{T_{\Ga g\Ga}}\arrow{s,l}
{\res^\Sigma_{\Sigma'}}
	\node{\bar H_q(\Sigma(g),V)}\arrow{s,r}{\res^{\Sigma(g)}_{\Sigma'(g)}}\\
\node{\bar H_q(\Sigma',V)}\arrow{e,t}{T_{\Ga g\Ga}}
	\node{\bar H_q(\Sigma'(g),V)}
\end{diagram}
$$
commutes.
\end{lemma}
\prf
Note first that, as $\Sigma$ is normal in $\Ga$, the group $\Ga$ acts on the finite set $\Ga g\Ga/\Sigma$ from the right by $\ga: h\Sigma\mapsto h\ga\Sigma$.
Therefore,
$$
|\Ga g\Ga/\Sigma|\=\sum_{h:\Ga g\Ga/\Ga}|\Ga_{h\Sigma}\bs\Ga|,
$$
where $\Ga_{h\Sigma}$ is the stabilizer of $h\Sigma$ in $\Ga$.
Note that 
$$
\ga\in\Ga_{h\Sigma}\ \Leftrightarrow\ h\ga\Sigma=h\Sigma\ \Leftrightarrow\ \ga\in\Sigma,
$$ 
so that
$$
|\Ga g\Ga/\Sigma|\=|\Ga g\Ga/\Ga||\Sigma\bs\Ga|.
$$
Let $v\in H_q^0(\Sigma, V)$.
Then $v$ represents the class $[v]=v+H_{q-1}^0(\Sigma,V)$ in $\bar H_q(\Sigma,V)$.
The element $\res^\Sigma_{\Sigma'}([v])$ is also represented by the same $v$, so that $T_{\Ga g\Ga}\res^\Sigma_{\Sigma'}([v])$ is represented by
$$
\frac 1{[\Ga:\Sigma']}\sum_{h:\Ga g\Ga/\Sigma'} hv.
$$
On the other hand, $T_{\Ga g\Ga}([v])$ is represented by
\begin{equation}
\frac 1{[\Ga:\Sigma]}\sum_i k_iv
\label{repre}
\end{equation}
where $\Ga g\Ga/=\dotcup_i k_i\Sigma$. Suppose that $\Sigma=\dotcup_j g_j\Sigma'$. Then, by Lemma 1.2.2(b), the sum
(\ref{repre}) equals
$$\frac 1{[\Ga:\Sigma]}\sum_i \left (
\frac 1{[\Sigma:\Sigma']}\sum_j k_i g_j v \right )
$$
modulo $H_{q-1}^0(\Sigma(g), V)$.
As $[\Ga:\Sigma']=[\Ga:\Sigma][\Sigma:\Sigma']$ and
$\Ga g\Ga=\dotcup_{i, j} k_i g_j\Sigma'$
the last sum equals $\frac 1{[\Gamma:\Sigma']}\sum_{\Ga g\Ga/\Sigma'} l v$
modulo $H_{q-1}^0(\Sigma(g),V)$ and its restriction to $\Sigma'$ equals the same sum modulo $H_{q-1}^0(\Sigma'(g),V)$.
\qed

In applications, it will be necessary to consider subsystems like the system of congruence subgroups defined as follows.
A subgroup $\Sigma$ of $\Ga$ is called a \e{congruence subgroup}, if it contains a subgroup of the form $\Ga(g_1)(g_2)\dots(g_n)$ for some $g_1,\dots g_n\in G$.

\begin{lemma}
The intersection of two congruence subgroups is a congruence subgroup.
\end{lemma}
\prf
We claim that
$$
\Ga(g_1)\dots(g_n)\ \bigcap\ \Ga(h_1)\dots(h_m)
$$
contains
$$
\Ga(g_1)\dots(g_n)(h_1)\dots(h_m).
$$
We spell out the argument in the case $m=n=1$ and leave the obvious iteration to the reader.
Recall that $\Sigma(g)$ is defined as the kernel of the map $\Sigma\to \Per(\Ga g\Ga/\Sigma)$.
So a given  $\ga\in\Ga$ lies in $\Ga(g)$ if and only if $\ga$ acts trivially on $\Ga g\Ga/\Ga$.
It lies in $\Ga(g)(h)$ iff it also acts trivially on $\Ga h\Ga/\Ga(g)$.
But then it acts trivially on $\Ga h\Ga/\Ga$, and hence $\ga$ lies in $\Ga(g)\cap\Ga(h)$.
\qed

Let
$$
L_q(\Ga,V)\=\varinjlim_\Sigma \bar H_q(\Sigma,V),
$$
where the limit is taken over all normal congruence subgroups $\Sigma$ of $\Ga$.
Further, let $L_q^{\rm all}(\Ga,V)$ denote the same direct limit, but now over all finite index normal  subgroups of $\Ga$.
The lemma above shows that one gets a well defined operator
$$
T_{\Ga g\Ga}: L_q(\Ga,V)\ \to\ L_q(\Ga,V), 
$$
and likewise for $L_q^{\rm all}(\Ga,V)$.
For the rest of the section, we consider the case $L_q(\Ga,V)$ only, but everything will as well apply to $L_q^{\rm all}(\Ga,V)$.

\begin{proposition}
The map $\1_{\Ga g\Ga}\mapsto T_{\Ga g\Ga}$ extends uniquely to a representation of the Hecke algebra $\CH(G,\Ga)$ on $L_q(\Ga,V)$.
\end{proposition}
\prf
Uniqueness is clear as the $1_{\Ga g\Ga}$ span the Hecke algebra.
We only have to show that the ensuing linear map is a representation.
For this we write it in a different manner.
Let $f\in\CH(G,\Ga)$ and let $v\in \bar H_q(\Sigma,V)$.
Define
$$
f.v\=\frac 1{[\Ga:\Sigma]}\sum_{y\in G/\Sigma}f(y) yv.
$$
Then $f.v$ lies in $\bar H_q(\Sigma',V)$, where $\Sigma'=\Sigma(g_1)\cap\dots\cap\Sigma(g_n)$ and the support of $f$ is contained in $\Ga g_1\Ga\cup\dots\cup\Ga g_n\Ga$.
This is an action of $\CH$ which extends the above map.
\qed

{\bf Remark.} A given $g\in G$ maps $H_q^0(\Sigma,V)$ to $H_q^0(g\Sigma g^{-1},V)\subset H_q^0(\Sigma\cap g\Sigma g^{-1},V)$ and thus it maps $\bar H_q(\Sigma,V)$ to $\bar H_q(\Sigma\cap g\Sigma g^{-1},V)$.
For $\Sigma'\subset\Sigma$ the diagram
$$\divide\dgARROWLENGTH by2
\begin{diagram}
\node{\bar H_q(\Sigma,V)}\arrow{e,t}g\arrow{s,r}{\res}
	\node{\bar H_q(\Sigma\cap g\Sigma g^{-1},V)}\arrow{s,r}{\res}\\
\node{\bar H_q(\Sigma',V)}\arrow{e,t}g
	\node{\bar H_q(\Sigma'\cap g\Sigma' g^{-1},V)}
\end{diagram}
$$
commutes.
Therefore the group $G$ acts on the limit $L_q(\Ga,V)$.
It is sometimes easier to understand $L_q(\Ga,V)$ as a $G$-module rather that a Hecke module.

\subsection{Unitary Hecke modules}
Suppose now that for every congruence subgroup $\Sigma$ the space $\bar H_q(\Sigma,V)$ is a Hilbert space in such a way that
\begin{itemize}
\item $\res^\Sigma_{\Sigma'}: \bar H_q(\Sigma,V)\to \bar H_q(\Sigma',V)$ is an isometry if $\Sigma'\subset \Sigma$, and
\item for each $g\in G$ the map $g:\bar H_q(\Sigma,V)\to\bar H_q(g\Sigma g^{-1},V)$, induced by the action of $g$ on $V$, is a unitary map.
\end{itemize}

The first condition gives $L_q(\Ga,V)$ the structure of a pre-Hilbert space.
The second implies that $G$ acts on this pre-Hilbert space by unitary maps.
If this is the case for every $q\ge 1$, we call $V$ a \e{unitary} Hecke module.

\begin{theorem}
Let $V$ be a unitary Hecke module.
For each $g\in G$ the operator $T_{\Ga g\Ga}$ is a bounded operator on the pre-Hilbert space $L_q(\Ga, V)$. The operator norm satisfies
$$
\norm T\ \le\ |\Ga g\Ga/\Ga|.
$$
The adjoint of $T_{\Ga g\Ga}$ is 
$$
T_{\Ga g\Ga}^*\=\frac{|\Sigma g\Sigma/\Sigma|}{|\Sigma g^{-1}\Sigma/\Sigma|}T_{\Ga g^{-1}\Ga}.
$$
In particular, if $\Ga g\Ga=\Ga g^{-1}\Ga$, then
$T_{\Ga g\Ga}$ is self-adjoint.

If there exists a locally compact group $G^*$ with $G\subset G^*$ such that $\Ga$ is a lattice in $G^*$, then we have
$$
T_{\Ga g\Ga}^*\=T_{\Ga g^{-1}\Ga}
$$
for every $g\in G$.
\end{theorem}
\prf
For $v\in \bar H_q(\Sigma,V)$ with $\norm v=1$ one has
\begin{eqnarray*}
\norm{T_{\Ga g\Ga} v} &=& \frac 1{[\Ga:\Sigma]}\norm{\sum_{h:\Ga g\Ga/\Sigma}hv}\\
&\le& \frac 1{[\Ga:\Sigma]}\sum_{h:\Ga g\Ga/\Sigma}\underbrace{\norm{hv}}_{\ \ \ \ =\norm v=1}\\
&=& \frac{|\Ga g\Ga/\Sigma|}{[\Ga:\Sigma]}\=|\Ga g\Ga/\Ga|.
\end{eqnarray*}
To compute the adjoint, let $v,w\in \bar H_q(\Sigma,V)$.
Then
\begin{eqnarray*}
\sp{T_{\Ga g\Ga}v,w}&=&\frac 1{[\Ga:\Sigma]}\sum_{h:\Ga g\Ga/\Sigma}\sp{hv,w}\\
&=&
\frac 1{[\Ga:\Sigma]}\sum_{h:\Sigma\bs\Ga g\Ga/\Sigma}\ \ \sum_{\sigma:\Sigma/\Sigma_{h\Sigma}}\sp{\sigma hv,w}\\
&=&
\frac 1{[\Ga:\Sigma]}\sum_{h:\Sigma\bs\Ga g\Ga/\Sigma}\ \ \sum_{\sigma:\Sigma/\Sigma_{h\Sigma}}\sp{hv,\sigma^{-1}w},
\end{eqnarray*}
where $\Sigma_{h\Sigma}$ is the stabilizer in $\Sigma$ of the coset $h\Sigma$ in the set $G/\Sigma$.
We have $\sigma^{-1} w\equiv w$ modulo $H_{q-1}^0$, so the last line equals
$$
\frac 1{[\Ga:\Sigma]}\sum_{h:\Sigma\bs\Ga g\Ga/\Sigma}|\Sigma h\Sigma/\Sigma|\sp{hv,w},
$$
where we have used $|\Sigma h\Sigma/\Sigma|=|\Sigma/\Sigma_{h\Sigma}|$.
Let $h=\ga_1 g\ga_2$ with $\ga_1,\ga_2\in\Ga$.
Then
\begin{eqnarray*}
|\Sigma h\Sigma/\Sigma| &=& |\Sigma \ga_1 g\ga_2\Sigma/\Sigma|\\
&=& |\ga_1\Sigma g\ga_2\Sigma/\Sigma|\\
&=& |\Sigma g\ga_2\Sigma/\Sigma|,
\end{eqnarray*}
as the map $x\Sigma\mapsto \ga_1 x\Sigma$ is a bijection on $G/\Sigma$.
Further, as $\Sigma$ is normal in $\Ga$, the group $\Ga$ acts on $G/\Sigma$ via $x\Sigma\mapsto x\ga\Sigma$.
Hence we get
$$
|\Sigma h\Sigma/\Sigma|\=|\Sigma g\Sigma/\Sigma|.
$$
And so
$$
\sp{T_{\Ga g\Ga}v,w}\=\frac{|\Sigma g\Sigma/\Sigma|}{[\Ga:\Sigma]}\sum_{h:\Sigma\bs \Ga g\Ga /\Sigma}\sp{v,h^{-1}w}.
$$
If $h$ runs through a set of representatives of $\Sigma\bs\Ga g\Ga/\Sigma$, then $h^{-1}$ runs through a set of representatives of $\Sigma\bs\Ga g^{-1}\Ga/\Sigma$.
It follows
$$
\sp{T_{\Ga g\Ga}v,w}\=\frac{|\Sigma g\Sigma/\Sigma|}{[\Ga:\Sigma]}\sum_{h:\Sigma\bs \Ga g^{-1}\Ga /\Sigma}\sp{v,hw}.
$$
Repeating the same argument with $g^{-1}$ in place of $g$ yields
$$
\sp{T_{\Ga g\Ga}v,w}\=\frac{|\Sigma g\Sigma/\Sigma|}{|\Sigma g^{-1}\Sigma/\Sigma|}\sp{v,T_{\Ga g^{-1}\Ga}w},
$$
or
$$
T_{\Ga g\Ga}^*\=\frac{|\Sigma g\Sigma/\Sigma|}{|\Sigma g^{-1}\Sigma/\Sigma|}T_{\Ga g^{-1}\Ga},
$$
 as claimed.
 
Suppose $\Ga g\Ga=\Ga g^{-1}\Ga$ and let $c=\frac{|\Sigma g\Sigma/\Sigma|}{|\Sigma g^{-1}\Sigma/\Sigma|}$.
We have to show that $c=1$.
With $T=T_{\Ga g\Ga}$ we have $T^*=cT$, so $T=(T^*)^*=cT^*=c^2T$.
As $c>0$ we conclude $c=1$.

Finally, assume the existence of $G^*$.
Note that $\Sigma$ is a lattice in $G^*$ as well and that
$$
|\Sigma g\Sigma/\Sigma|\=|\Sigma/\Sigma\cap g\Sigma g^{-1}|,
$$
as the set $\Sigma g\Sigma/\Sigma$ is one orbit under the left translation action of $\Sigma$ and the stabilizer of the point $g\Sigma/\Sigma$ is $\Sigma\cap g\Sigma g^{-1}$.
Accordingly,
$$
|\Sigma g^{-1}\Sigma/\Sigma|\=|\Sigma/\Sigma\cap g^{-1}\Sigma g|\=|g\Sigma g^{-1}/\Sigma\cap g\Sigma g^{-1}|.
$$
Let $\mu$ be the left Haar measure on $G^*$, then $\mu(g\Sigma g^{-1}\bs G^*)=\mu(\Sigma\bs G^*)$ and so
{ \begin{multline*}
\mu\(G^*/(\Sigma\cap g\Sigma g^{-1})\) 
|g\Sigma g^{-1}/\Sigma\cap g\Sigma g^{-1}|\=\mu(G^*/g\Sigma g^{-1})\=\mu(\Sigma\bs G^*)\\
\=\mu\(G^*/(\Sigma\cap g\Sigma g^{-1})\)|\Sigma /\Sigma\cap g\Sigma g^{-1}|.
\end{multline*}}
This implies the claim.
\qed

A Hecke pair $(G,\Ga)$ is called \e{unimodular}, if 
$$
|\Ga g\Ga/\Ga|\= |\Ga\bs\Ga g\Ga|
$$
holds for every $g\in G$.

\begin{corollary}
Let $(G,\Ga)$ be a Hecke pair such that there exists a locally compact group $G^*$ as in the theorem, then the pair $(G,\Ga)$ is unimodular.
\end{corollary}

\prf
The last part of the proof of the theorem gives $|\Ga g\Ga/\Ga|=|\Ga g^{-1}\Ga/\Ga|$.
The inversion $h\mapsto h^{-1}$ induces a bijection between $\Ga g^{-1}\Ga/\Ga$ and $\Ga\bs\Ga g\Ga$.
\qed

Note that the condition
$$
\Ga g\Ga\=\Ga g^{-1}\Ga
$$ is satisfied in important examples like the Hecke pairs $(\SL_2(\Q),\SL_2(\Z))$ or $(\PGL_2(\Q),\PGL_2(\Z))$.
Note further, that this condition implies that the Hecke algebra $\CH(G,\Ga)$ is commutative, i.e., that $(G,\Ga)$ is a \e{Gelfand pair}.

\begin{example}
For the sake of completeness we give an example of a Hecke pair $(G,\Ga)$ which is not unimodular.
Let $p$ be a prime and let $G$ be the semidirect product $\Q_p\rtimes\Q_p^\times$.
So $G$ is the topological space $\Q_p\times \Q_p^\times$ with the multiplication $(x,y)(x',y')=(x+yx',yy')$.
Let $\Ga$ be the compact open subgroup $\Z_p\rtimes\Z_p^\times$ and let $g=(0,p)$.
Then $|\Ga g\Ga/\Ga|=p$, whereas $|\Ga\bs\Ga g\Ga|=1$.
\end{example}

{\bf Remark.} Note that if $V$ is a unitary Hecke module, then the representation of $G$ on the pre-Hilbert space $L_q(\Ga,V)$ is unitary.

\subsection{Lowering the order}
There is a canonical injective linear map
$$
l_q=l_{V,\Sigma,q}: \bar H_q(\Sigma,V)\ \hookrightarrow\ \Hom(\Ga,\C)\otimes \bar H_{q-1}(\Sigma,V)
$$
given as follows: first note that there is a canonical isomorphism $\Hom(\Ga,\C)\otimes \bar H_{q-1}(\Sigma,V)\tto\cong \Hom(\Ga, \bar H_{q-1}(\Sigma,V))$.
Using this, we can define $l_q$ as
$$
l_q(v)(\ga)\= (\ga-1)v.
$$
This indeed is well defined as for $\ga,\tau\in\Sigma$ one has $(\ga\tau-1)\equiv (\ga-1)+(\tau-1)\mod I_\Sigma^2$.
We call $l_q$ the \e{order lowering homomorphism}.
It can be used to establish a unitary structure as is shown in the next section.

In the following sections we will give examples for spaces $V$ to which the Hecke calculus applies.
For each of these spaces the following problems arise.
\begin{itemize}
\item Determine the spectral decomposition of $T_{\Ga g\Ga}$ on the Hilbert completion of $L_q(\Ga,V)$.
If the Hecke algebra is commutative, one can give a simultaneous spectral decomposition.
Are the Hecke operators compact?
\item Determine the difference that it makes for the spectrum of the Hecke algebra whether one starts with congruence subgroups only or all finite index normal subgroups $\Sigma$.
This problem should be related to the Selberg conjecture.
\item Determine the action of the Hecke algebra in terms of Fourier-Taylor expansions (defined below) and on the ensuing $L$-functions.
\end{itemize}

Instead of the Hecke-action, it is sometimes more useful to consider the $G$-action on $L_q(\Ga,V)$.
Note that the order-lowering homomorphism can be iterated to 
$$
\bar H_q(\Sigma,V)\ \hookrightarrow\ \Hom(\Sigma,\C)^{\otimes q}\otimes H^0(\Sigma, V).
$$
Let 
$$
H_\Ga\df \lim_{\stack\to\Sigma}\Hom(\Sigma,\C).
$$
Then $H_\Ga$ is a $G$-module and, taking limits,  one gets an injection of $G$-modules,
$$
L_q(\Ga,V)\ \hookrightarrow\ (H_\Ga)^{\otimes q}\otimes V^\infty,
$$
where $V^\infty=\bigcup_{\Sigma} V^\Sigma$.
Sometimes this map will be surjective on the Hilbert space completions.
Thus the spectral problem decomposes into
\begin{itemize}
\item determining the spectral decomposition of the module $H_\Ga$, resp. its Hilbert space completion, and
\item determining the decomposition of tensor products.
\end{itemize}

These aims will be pursued in future work.

\section{Automorphic forms}

\subsection{Modular forms}
For any ring $R$, let $G_R=\SL_2(R)$.
The Lie group $G_\R$ acts by linear fractionals on the upper half plane $\H$ in $\C$.
Let $\Ga\subset G_\R$ be a lattice.
A \e{cusp} of $\Ga$ is a parabolic fixed point in the boundary $\partial\H=\hat\R\=\R\cup\{ \infty\}$.
For a cusp $c$, let $\Ga_c$ be its stabilizer in $\Ga$.
For every cusp $c$ there exists an element $\sigma_c^\Ga\in G_\R$ such that $\sigma_c^\Ga\infty=c$ and the action of $(\sigma_c^\Ga)^{-1}\Ga_c\sigma_c^\Ga$ on $\H$ is generated by the translation $z\mapsto z+1$.
For $k\in 2\Z$ and $f:\H\to\C$ define
$$
(f|_k\ga)(z)\= (cz+d)^{-k}f(\ga z),
$$
where $\ga=\left(\begin{array}{cc}* & * \\c & d\end{array}\right)\in G_\R$.
Since we rather deal with left actions of the groups and we keep $k$ fixed, we also write $\ga f=f|_k\ga^{-1}$.

Let $\Ga=\SL_2(\Z)$, then $(G_\Q,\Ga)$ is a Gelfand pair.
Let $V=\CM_k$ be the complex vector space of all holomorphic functions $f:\H\to\C$ such that for every cusp $c$ of $\Ga$ the function $f|_k\sigma_c^\Ga$ is bounded on the set $\{ z\in\H:\Im(z)\ge 1\}$.

A congruence group $\Sigma\subset \Ga$ is a subgroup which contains the group
$$
\Ga(N)\=\ker(\Ga\to\SL_2(\Z/N\Z))
$$
for some $N\in\N$.
Note that this definition coincides with our earlier definition of a congruence group for the pair $(G_\Q,\Ga)$.
Let $\Sigma\subset\Ga$ be a congruence subgroup.
The space $\CM_k(\Sigma)=V^\Sigma=H^0(\Sigma,V)$ is the space of holomorphic modular forms for the group $\Sigma$.
For every $f\in\CM_k(\Sigma)$ and every cusp $c$ of $\Ga$ the function $z\mapsto f|_k\sigma_c^\Ga(z)$ is periodic of period $|\Ga_c/\Sigma_c|$.
It therefore has a Fourier-expansion.
The normalized zeroth coefficient
$$
d_0(f,c)\=\frac 1{|\Ga_c/\Sigma_c|}\int_0^{|\Ga_c/\Sigma_c|} f|_k\sigma_c^\Ga(z+t)\, dt
$$
is independent of $z\in\H$.
Let $S_k(\Sigma)$ be the subspace of all cusp forms $f$ for $\Sigma$, 
which are characterized by the fact that for every cusp $c$ of $\Ga$  one has
$$
d_0(f,c)\= 0.
$$
Note that every cusp form $f$ has rapid decay at the cusps, which means 
that for every cusp $c$ one has $f|_k\sigma_c^\Ga(x+iy)=O(e^{-\al y})$ 
for some $\al>0$.

\begin{theorem}
There is a canonical choice of inner products that makes $V$ a unitary Hecke module.
\end{theorem}
\prf
As a first step we have to construct a canonical inner product on $\CM_k(\Sigma)$.
For this note that for $f\in\CM_k(\Sigma)$ and $g\in S_k(\Sigma)$ one has the Petersson inner product:
$$
\sp{f,g}_{\rm Pet}\=\frac 1{\vol{(\Sigma\bs\H)}}\int_{\Sigma\bs\H}f(z)\ol{g(z)}y^k\,d\mu(z),
$$
where $\mu$ is the $G_\R$-invariant measure $\frac{dxdy}{y^2}$.
We extend it to an inner product on $\CM_k(\Sigma)$ as follows.
The Petersson inner product defines an orthogional projection $P:\CM_k(\Sigma)\to S_k(\Sigma)$.
Let $C$ be the set of cusps of $\Ga$.
For $f,g\in \CM_k(\Sigma)$ let 
$$\sp{f,g}_{\rm cusp}\=\frac 1{|\Sigma\bs C|}\ \sum_{c\in \Sigma\bs C}d_0(f,c)\ol{d_0(g,c)},
$$
where the sum runs over all $\Sigma$-equivalence classes of cusps $c$.
The form
$$
\sp{f,g}\=\sp{f,g}_{\rm cusp}+\sp{P(f),P(g)}_{\rm Pet}
$$
is a positive definite inner product such that for $\Sigma'\subset\Sigma$ the inclusion $\CM_k(\Sigma)\hookrightarrow\CM_k(\Sigma')$ is an isometry and for $g\in G_\Q$ the map $g:\CM_k(\Sigma)\to\CM_k(g\Sigma g^{-1})$ is unitary.
This settles the case $q=1$.
Suppose the case $q$ has been taken care of, then consider the order lowering homomorphism
$l_{q+1}:\bar H_{q+1}(\Sigma,V)\hookrightarrow \Hom(\Sigma,\C)\otimes \bar H_q(\Sigma,V)$.
The Eichler-Shimura map gives a canonical isomorphism
$$
\Hom(\Sigma,\C)\ \cong\ \CM_2(\Sigma)\oplus\ol{S_2(\Sigma)}.
$$
The right hand side is equipped with a canonical inner product by the above,
such that for $\Sigma'\subset \Sigma$ one has a commutative diagram
$$
\begin{diagram}
\node{\Hom(\Sigma,\C)}\arrow{e,J}\arrow{s,J}
	\node{\CM_2(\Sigma)\oplus\ol{S_2}(\Sigma)}
		\arrow{s,J}\\
\node{\Hom(\Sigma',\C)}\arrow{e,J}
	\node{\CM_2(\Sigma')\oplus\ol{S_2}(\Sigma').}
\end{diagram}
$$
The Petersson inner product thus establishes a canonical inner product on $\Hom(\Sigma,\C)$ such that the restriction map $\Hom(\Sigma,\C)\to\Hom(\Sigma',\C)$ is isometric and for $g\in G_\Q$ the map $g:\Hom(\Sigma,\C)\to\Hom(g\Sigma g^{-1},\C)$ is unitary.
Now equip $\bar H_{q+1}(\Sigma,\CM_k)$ with the inner product from the above injection into $\Hom(\Sigma,\C)\otimes \bar H_q(\Sigma,V)$ to make $\CM_k$ a unitary Hecke module in a canonical way.
\qed

The inner product $<\cdot, \cdot>$ on  $\bar H_q(\Sigma,\CM_k)$
can inductively be written explicitly as follows.
We write $\tilde M_k^q=\tilde M_k^q(\Sigma)$ for the space $H_q^0(\Sigma,\CM_k)$.
Then $\bar H_q(\Sigma, \CM_k)=\tilde M_k^q/\tilde M_k^{q-1}$.
 Assume that the inner product $<\cdot, \cdot>$ has been defined on
$\tilde M_k^{q-1}/\tilde M_k^{q-2}$  and let $\{f_i\}_{i=1}^d$ be an
orthonormal basis of $\tilde M_k^{q-1}/\tilde M_k^{q-2}$.

Let now $f, g \in \tilde M_k^{q}/ \tilde M_k^{q-1}$ and suppose that for every $\sigma \in \Sigma$,
$$
f|_k(\sigma-1)=\sum_{i=1}^d \phi_i(\sigma) f_i
$$
$$
g|_k(\sigma-1)=\sum_{i=1}^d \psi_i(\sigma) f_i
$$
for some $\phi_i, \psi_i \in$ Hom$(\Sigma, \C).$
Suppose that the Eichler-Shimura isomorphism maps $\phi_i$ (resp. $\psi_i$) to $(a_i, \bar b_i)$ (resp. $(c_i, \bar d_i)$) with $a_i, 
c_i \in \CM_2$ and $b_i, d_i \in S_2$. Then, the definition of inner
product on $\bar H_{q}(\Sigma,\CM_k)$
can be written as
\begin{eqnarray*}
<f, g> &=& \sum_{i, j}  \Big (<a_i, c_j><f_i, f_j>+\overline{<b_i,
d_j>}<f_i, f_j> \Big ) \\
&=&\sum_i  \Big (<a_i, c_i>+\overline{<b_i, d_i>}
\Big
).
\end{eqnarray*}
It should be noted that \"O. Imamoglu and C. O'Sullivan \cite{IO} have 
given an
alternative definition of an inner product on a subspace of $\tilde 
M_k^1/\tilde M_k^0$
and (conjecturally) for $\tilde M_k^q/\tilde M_k^{q-1},$ ($q\ge 2$). 
Their construction relies on delicate analytic manipulations but the inner 
product it gives can be proved to have essentially the same value as our 
inner product.

\subsection{Cusp forms}
In this section, let $\Ga$ be a congruence subgroup of $G_\Z$.

By linearity, we extend the definition $f|_k\sigma$ to elements $\sigma$ of
the group ring $\R[\Ga]$.
Let $k\ge 0$ be even and let $S_k(\Ga)$ be the space of cusp forms of
weight $k$.

We now define cusp forms of order $q$, which are the cuspidal analogues of the elements of 
$\tilde M_k^q$ discussed at the end of Sec. 2.1.
First let $\tilde S_{k,0}(\Ga)=S_{k,0}(\Ga)=S_k(\Ga)$, so classical cusp forms are of order $0$.
Next suppose $\tilde S_{k,q}(\Ga)$ and $S_{k,q}(\Ga)$ are already 
defined and let $\tilde S_{k,{q+1}}(\Ga)$  be the space of all functions $f$ with
\begin{itemize}
\item $f:\H\to\C$ holomorphic,
\item $f|_k(\ga -1)\in \tilde S_{k,q}(\Ga)$ for every $\ga\in\Ga$,
\item for every cusp $c$, $(f|_k\sigma_c^\Ga)(z)=O(e^{-\al y})$ as $y\to\infty$
for some $\al >0$ (``rapid decay at the cusps").
\end{itemize}
Further, let $S_{k,{q+1}}(\Ga)$ be the set of all $f\in \tilde S_{k,{q+1}}(\Ga)$ with
$f|_k(\ga-1)=0$ for every parabolic element $\ga$ of $\Ga$.
Note that for $f\in S_{k,{q+1}}(\Ga)$ and $\ga\in\Ga$ one has $f|_k(\ga-1)\in S_{k,q}(\Ga)$.

Compare $\tilde S_{k,1}$ with $PS_{k,2}$ of \cite{DKMO} where a 
classification and a converse-theorem-type of result for such 
functions is proved.

Note that the space $V$ of all holomorphic functions on $\H$ which satisfy the above growth condition at each cusp, serves as a 
$G$-module, where, for instance, $G=\SL_2(\Q)$ and $\Ga=\SL_2(\Z)$.
In this way one obtains a theory of Hecke operators for the space $\tilde S_{k,q}(\Ga)/\tilde S_{k,{q-1}}(\Ga)$. Also, in view of
the discussion at the end of Section 2.1., we have an inner product on the same space that makes it a unitary Hecke module.

\subsection{Fourier-Taylor expansion}
Set $T=\left(\begin{array}{cc} 1 & 1 \\ 0 & 1\end{array}\right)$ and let 
$V_0$
be the space of all holomorphic functions $f$ on the upper half plane with
$f_0|T=f$, i.e. periodic with period $1$.
Inductively, let $V_{q+1}$ be the space of all holomorphic functions on
$\H$ such that $f|_0(T-1)\in V^q$.
We also set $V_{-1}=\{0\}$.
Note that $S_{k,q}(\Ga)\subset V_q$ if $\Ga$ contains the translation $z\mapsto z+1$.

\begin{proposition}
Every $f\in V_q$ has a Fourier-Taylor-expansion
$$
f(z)\=\sum_{n\in\Z}e^{2\pi inz}(a_{n,0}+a_{n,1}z+\dots+a_{n,q}z^{q})
$$
for uniquely determined coefficients $a_{k,j}\in\C$.
For every $j$ and every $y>0$ the sequence $(a_{n,j}e^{-2\pi ny})_{n\in\N}$ is rapidly decreasing.
The map $T:f\mapsto f(z+1)-f(z)$ is a surjection from $V_q$ to $V_{q-1}$.
With the natural inclusion $V_0\hookrightarrow V_q$ one gets an exact sequence
$$
0\to V_0\to V_q\to V_{q-1}\to 0.
$$
\label{Fourtayl}
\end{proposition}
\prf
By induction on $q$.
For $q=0$, every $f\in V_0$ is periodic, therefore has a Fourier-expansion, which, as $f$ is
holomorphic, is of the form $f(z)=\sum_na_ne^{2\pi inz}$.
Next assume the claim proven for $q$ and let $f\in V_{q+1}$.
By induction,
$$
f(z+1)-f(z)\=\sum_n e^{2\pi inz}(a_{n,0}+\dots+a_{n,q}z^{q}).
$$
There is a unique sequence of numbers $b_{n,j}$ that satisfy the recursion relations
$$
a_{n,k}\=\sum_{j=k+1}^{q+1} \left(\begin{array}{c}j \\k\end{array}\right) b_{n,j}
$$
for $k=0,\dots,q$. For each $y>0$ the sequence $b_{n,j}e^{-2\pi ny}$ is rapidly decreasing.
Let $g(z)=\sum_ne^{2\pi inz}(b_{n,1}z+\dots+b_{n,q+1}z^{q+1})$,
then $g$ is holomorphic, and
$$
g(z+1)-g(z)\=\sum_n e^{2\pi inz}(a_{n,0}+\dots+a_{n,q}z^{q})\= f(z+1)-f(z).
$$
So $h(z)=f(z)-g(z)\in V_0$.
As $g$ and $h$ possess Fourier-Taylor expansions, so does $f$.
By induction, the coefficients $a_{n,k}$ are uniquely determined,
so are those of $h(z)$, which implies that the coefficients of $f$ are uniquely determined.
\qed

\subsection{Intervention of Lie groups}
For $\left(\begin{array}{cc}a & -b \\b & a\end{array}\right)\in K_\R=\SO(2)$ define
$$
\eps_k\left(\begin{array}{cc}a & -b \\b & a\end{array}\right)\=(a+ib)^{-k}.
$$
Let $G_\R=ANK_\R$ be the Iwasawa decomposition and let $\ul k:G_\R\to K_\R$ be the corresponding projection. Then
$$
\ul k\left(\begin{array}{cc}a & b \\c & d\end{array}\right)\= \frac 1{\sqrt{c^2+d^2}}
\left(\begin{array}{cc}d & -c \\c & d\end{array}\right).
$$
For a given function $f$ on the upper half plane define the function $\psi_f$ on $G_\R$ by
$$
\psi_f(g)=(\Im (gi))^{k/2}\,\eps_{k}(\ul k(g))\,f(gi).
$$
The next lemma shows that, via the identification $\psi$, the action
$|_k$ on functions on $\H$ becomes action by left translation.
\begin{lemma}
For $\ga\in G_\R$ and $x\in G_\R$ one has
$$
\psi_f(\ga x)\=\psi_{f|_k\ga}(x).
$$
\label{corr}
\end{lemma}
\prf A computation relying on the identity
$$
\eps_k(\ul k(\ga g))\=\(\frac{cz+d}{| cz+d|}\)^{-k}\eps_k(\ul k(g)),
$$
where $z=gi\in\H$.
\qed

We next define a sequence of Hilbert spaces $L^{2}_{q}(\Ga\bs G)$, where $q=0,1,2,\dots$.
We start with providing spaces of measurable functions.
Define $\CF_{-1}(\Ga\bs G)=0$ and let $\CF_0(\Ga\bs G)$ be the space of all $\Ga$-invariant measurable functions on $G$.
Inductively, define $\CF_{q+1}$ to be the set of all measurable functions $f$ on $G$ with $(\ga -1)f\in\CF_q$ for every $\ga\in\Ga$.
Let $\bar\CF_{q+1}=\CF_{q+1}/\CF_q$.
Consider the map $\eta:\CF_{q+1}\to\Hom(\Ga,\bar\CF_q)$ given by
$$
\eta(f)(\ga)\=(\ga-1)f.
$$
This indeed defines a group homomorphism in $\ga$, as for $\ga,\tau\in \Ga$ one has $(\ga\tau-1)\equiv(\ga-1)+(\tau-1)$ modulo $I^2$, and $I^2\CF_{q+1}\subset\CF_{q-1}$.
The kernel of $\eta$ equals $\CF_q$, so $\eta$ defines an injection $\bar\CF_{q+1}\hookrightarrow \Hom(\Ga,\bar\CF_q)\cong\Hom(\Ga,\C)\otimes\bar\CF_q$.
By iteration one gets an injective map
$$
\eta_{q}:\bar\CF_{q}\ \hookrightarrow\ \Hom(\Ga,\C)^{\otimes q}\otimes\CF_0,
$$
where we have used $\bar\CF_0=\CF_0$.

\begin{lemma}
Suppose that $\Ga$ is torsion-free.
Then $\eta_{q}$ is surjective.
\end{lemma}

\prf
If the group $\Ga$ has genus $g$ and $s$ cusps, then there are
$2g$ hyperbolic elements $\gamma_i$ and $s$ parabolic elements
$\ga_{2g+i}$ generating $\Gamma$ and satisfying the relation:
$$[\gamma_1, \gamma_{g+1}]\dots
[\gamma_{g}, \gamma_{2g}] \ga_{2g+1}\cdots\ga_{2g+s}=1.$$

Because of this relation, every element in $\Hom(\Ga,\C)^{\otimes q}\otimes\CF_0$ is uniquely determined by its values at $(\ga_{i_1}, \dots
\ga_{i_q})$ for $i\in \{1,\dots,N\}^q$, where $N=\max(2g,2g+s-1)$.
Therefore, to
establish surjectivity, it suffices to show that for every choice of 
functions $f_{i_1,
\dots, i_q} \in \CF_0$ there is an $F \in \bar\CF_{q}$ such that
\begin{eqnarray}(\ga_{i_q}-1) \dots
(\ga_{i_1}-1)F=f_{i_1, \dots, i_q}.\label{surj}
\end{eqnarray}
By proposition 4.1 of \cite{DSr}, for every $q$-tuple of integers
$L=(l_1, \dots, l_q)$ in $\{1, \dots, N\}$ there is a smooth
function $\Lambda_L$ on $\H$ such that $$\Lambda_{L}|_0(\gamma_{i_1}-1)\dots
(\gamma_{i_q}-1)=(-1)^{q} \delta_I^L$$
for any $s$-tuple $I$, where $\delta_I^L$ is the Kronecker delta
function of the $q$-tuple, namely $\prod_k \delta_{i_k}^{l_k}$.
Then, with Lemma \ref{corr}, the function
$$F:=\sum_{L} f_{L} \psi_{\Lambda_L}$$
with the sum ranging over all $L$ with components in $\{1, \dots, N\}$
satisfies (\ref{surj}).
\qed

Next, let $\CN_0$ denote the space of all nullfunctions on $G$.
Let
$$
\CN_{q}\df\eta_{q}^{-1}\(\Hom(\Ga,\C)^{\otimes q}\otimes\CN_0\).
$$
Let $\CL_{q}^2\= \eta_{q}^{-1}\(\Hom(\Ga,\C)^{\otimes q}\otimes\CL^2(\Ga\bs G)\)$, where $\CL^2(\Ga\bs G)$ is the space of all square integrable functions on $\Ga\bs G$.
Finally, we define
$$
L^2_{q}(\Ga\bs G)\df \CL_{q}^2/\CN_{q}.
$$
The space
$L^2(\Ga\bs G)=\CL^2(\Ga\bs G)/\CN_0$ is a Hilbert space.
Fix a Hilbert space structure on $\Hom(\Ga,\C)$ and equip the space $L_{q}^2(\Ga\bs G)$, which is mapped bijectively onto $\Hom(\Ga,\C)^{\otimes q}\otimes L^2(\Ga\bs G)$ the induced Hilbert space structure.

\subsection{No intervention of ad\`eles}
Let $\A=\A_\fin\times\R$ be the adele-ring over $\Q$.
Let $K_\Ga$ be a compact open subgroup of $G_{\A_\fin}$, then $\Ga=K_\Ga\cap G_\Q$
is a congruence subgroup, and the natural map
$$
\Ga\bs G_\R\to G_\Q\bs G_\A /K_\Ga,
$$
which maps $\Ga x$ to $G_\Q x K_\Ga$, is a $G_\R$-equivariant, continuous bijection.
This gives a natural isomorphism
$$
L^2(\Ga\bs G)\ \to\ L^2(G_\Q\bs G_\A)^{K_\Ga}.
$$
In other words, automorphic forms on $\Ga\bs G$ can be lifted to $G_\Q\bs G_\A$.
This is what Dieudonn\'e calls the ``intervention of adeles''.
We ask for higher forms in the adelic setting.
Note first that there are no higher $G_\Q$-invariants, as
 the group $G_\Q$ is perfect.

Also the $K_\Ga$-action does not yield higher order forms, at least not in the space $L^2(G_\Q\bs G_\A)$, as the group $K_\Ga$ is 
compact and acts through a continuous representation on the Hilbert space $L^2(G_\Q\bs G_\A)$, see Proposition \ref{1.1}.

\section{L-functions}
\subsection{Higher order cusp forms}
Let $\Ga$ be a lattice in $G$ such that $\infty$ is a cusp of $\Ga$ of 
width 1.
Let $f\in\tilde S_{k,q}(\Ga)$.
Then $f$ has the Fourier-Taylor expansion
$$
f(z)\=\sum_{n=1}^{\infty}e^{2\pi in {z}}(a_{n,0}+a_{n,1} 
z+\dots+a_{n,q}z^{q}).
$$
For $\nu=0,\dots,q$ define
$$
L_\nu(f,s)\=\sum_{n=1}^\infty a_{n,\nu}\,n^{-s}.
$$
Let
\begin{eqnarray*}
\La(f,s) &=& \int_0^\infty f(iy)\, y^{s-1}\,dy\\
&=& \sum_{\nu=0}^{q} i^\nu \Ga(s+\nu)\({2\pi}\)^{-(s+\nu)} L_\nu(f,s+\nu).
\end{eqnarray*}
Let $S=\left(\begin{array}{cc}0 & -1 \\1 & 0\end{array}\right)$.
After replacing $\Ga$ with a conjugate if necessary, we may assume that $0$ is a cusp of $\Ga$, too.
Then the group $S^{-1}\Ga S$ has $\infty$ for a cusp.
Let $w=w_{\Ga}>0$ be its width.
Let $S_{w}$ be the matrix $S\diag(\sqrt{w},\sqrt{w}^{-1})$.
Then $\infty$ is a cusp of width $1$ of the group
$$
\hat\Ga\= S_{w}^{-1}\Ga S_{w}.
$$
Let $\hat f=f|_k S_{w}$.
Then $\hat f\in \tilde S_{k,q}(\hat\Ga)$.
Note that $S_{w}^{2}=1$, and so $w_{\hat\Ga}=w_\Ga$ and $\hat{\hat\Ga}=\Ga$ as well as $\hat{\hat f}=f$.

As an example let 
$N\in\N$ and consider the group $\Ga=\Ga_0(N)$ consisting of all 
$\left(\begin{array}{cc}a & b \\c & d\end{array}\right)\in\SL_2(\Z)$ such 
that $c\equiv 0\mod N$. In this case one has $w=N$ and 
$\hat\Gamma=\Gamma$.

\begin{proposition}\label{1.2}
Assume that $0$ and $\infty$ are cusps of $\Ga$.
For $f\in\tilde S_{k,q}(\Ga)$ the function $\La(f,s)$ extends to an entire function and satisfies the functional equation
$$
\La(f,s)=i^k w^{s-\frac k2}\La(\hat f,k-s).
$$
For
$0\le\nu\le{q-1}$ the $L$-function $L_\nu(f,s)$ extends to an entire function.
\end{proposition}
\prf
Note that $\hat f(z)=(\sqrt{w}z)^{-k}f(\frac{-1}{wz})$, 
so in particular one has $f(i\frac 1{wy})=(i\sqrt w y)^k\hat f(iy)$.
We decompose the integral defining $\La(f,s)$ as $\int_0^{1/\sqrt w}+\int_{1/\sqrt w}^\infty$.
The second gives an entire function $F(s)$.
The first is
\begin{eqnarray*}
\int_0^{1/\sqrt w}f(iy)\,y^{s}\,\frac{dy}y &=& w^{-s} \int_{1/\sqrt w}^\infty f\(i\frac 1{wy}\)y^{-s}\,\frac{dy}y\\
&=& 
i^kw^{-s+\frac k2}\int_{1/\sqrt w}^\infty \hat f(iy)y^{k-s}\,\frac{dy}y\\
&=& i^k w^{-s+\frac k2}\hat F(k-s),
\end{eqnarray*}
where $\hat F$ is the same as $F$ with $\hat f$ in place of $f$.
Therefore, $\La(f,s)=F(s)+i^k w^{\frac k2-s}\hat F(k-s)$, which proves that $\La(f,s)$ is entire and satisfies the functional quation.

Note that, with $\Lambda_j(f, s):=w^s\Gamma(s)(2 \pi)^{-s} L_j(f, s),$
we have that
$
\La(f,s)$ equals $\sum_{\nu=0}^{q} i^\nu\La_\nu(f,s+\nu),
$
Form which the second claim follows in an easy induction.
\qed

\subsection{Convolution of L-functions}
In \cite{CDO}, it is shown that if
$f(z)=\sum_{n=1}^\infty a_ne^{2\pi i nz}$ and
$g(z)=\sum_{n=1}^\infty b_ne^{2\pi i nz}$
are cusp forms of weight $k \in \mathbb Z_{\ge 2}$ and $2$ respectively,
then the $n$-th Fourier coefficient of the second order form
$
F(z)=f(z) \int_{i \infty}^z g(w)dw
$ is
$$\sum_{j=1}^{n-1} \frac{a_{n-j}b_j}{j}.$$
Following \cite{DKMO}, one can define the $L$-function of $F(z)$ by means of
the Dirichlet series
$$\sum_{n=1}^\infty n^{-s}\sum_{j=1}^{n-1} \frac{a_{n-j}b_j}{j}.$$
This function of $s$ admits meromorphic
continuation and satisfies a functional equation (\cite{DKMO}).

Again we assume that $0$ and $\infty$ are cusps of $\Ga$, the width of $\infty$ being $1$.
We define $w$, and $\hat\Ga$ as in Theorem \ref{1.2}.
Let $k,l$ be even integers $\ge 0$ and let $f \in S_k(\Ga), g \in S_l(\Ga).$ 
Denote their respective Fourier coefficients by $a_n$ and $b_n$ respectively.
For complex numbers  $s$ and $t$ with large enough real parts we set
$$
(L_f\# L_g)(s,t)\=\sum_{n=1}^\infty n^{-s}\sum_{j=1}^{n-1}
\frac{a_{n-j}b_j}{j^t}.
$$
For Re$(s)$ large enough, we observe that the Mellin transform
\begin{eqnarray*}
\La_{f,g}(s,t) &=& \int_{0}^{\infty} f(ix)\int_{0}^{\infty}
g(ix+iy)y^{t-1}dy\,x^{s-1}dx
\end{eqnarray*}
equals
\begin{multline}
 \sum_{n, m} a_n b_m\int_{0}^{\infty} \int_{0}^{\infty}
e^{-2\pi n x} e^{-2\pi m x-2\pi m y}y^{t-1}x^{s-1}dxdy=\\
(\frac{1}{2\pi})^{s+t} \Gamma(s)\Gamma(t) \sum_{n, m}
\frac{a_nb_m}{(n+m)^sm^t}=
(\frac{1}{2\pi})^{s+t}\Gamma(s)\Gamma(t) (L_f \# L_g)(s,t)
\end{multline}
We shall analytically continue $(L_f \# L_g)(s,t)$  by a repeated application of the Riemann trick.
We decompose the integrals into sums of the form $\int_0^{1/\sqrt w}+\int_{1/\sqrt w}^\infty$.
Then $\La_{f,g}(s,t)$ will be the sum of four terms,
$$
\underbrace{\int_0^{1/\sqrt w}\int_0^{1/\sqrt w}}_{=A}+\underbrace{\int_0^{1/\sqrt w}\int_{1/\sqrt w}^\infty}_{=B}+\underbrace{\int_{1/\sqrt w}^\infty\int_0^{1/\sqrt w}}_{=C}+\int_{1/\sqrt w}^\infty\int_{1/\sqrt w}^\infty.
$$
The last summand defines a holomorphic function on $\C^2$, as $f$ and $g$ are rapidly decreasing at $\infty$.

Recall that we have $f(i\frac 1{wy})= (i \sqrt w y)^k\hat f(y)$ and 
likewise for $g$ and $l$.
The substitution $x\mapsto 1/wx$ in the outer integral shows that $A$ equals
\begin{multline*}
w^{-s}\int_{1/\sqrt w}^\infty\int_0^{1/\sqrt w}f(i\frac 
1{wx})g \Big (i(\frac 1{wx}+y) \Big ) y^{t-1}x^{-s-1}\,dy\,dx=\\
i^kw^{-s+k/2}\int_{1/\sqrt w}^\infty\int_0^{1/\sqrt w}\hat f(ix)g \Big (i(\frac 
1{wx}+y) \Big ) y^{t-1}x^{k-s-1}\,dy\,dx
\end{multline*}
As $f$ and $g$ are rapidly decreasing at $\infty$, this integral converges absolutely for $\Re(t)>0$ and defines a holomorphic function in $(s,t)$ in that region.
The inner integral can be integrated by parts to get
$$
\frac 1t g\(i(\frac 1{wx}+y)\)y^t|_0^{1/\sqrt w}-\frac 1t\int_0^{1/\sqrt w}y^t\frac{\partial}{\partial y}g\(i(\frac 1{wx}+y)\)\,dy.
$$
So $A$ extends to a meromorphic function on $\Re(t)>-1$ with a simple pole at $t=0$.
Iteration of that argument shows that $A$ is entire except for simple poles at $t=0,-1,-2,\dots$.

The other summands can be treated similarly.
Since the $\Ga$-function has poles  at $t=0,-1,-2,\dots,$ we have proven the following theorem.

\begin{theorem}\label{ac}
The function $L_f\# L_g$ extends to an entire function on $\C^2$.
\end{theorem}

If $t=p\in\N$ we set
$(L_f \#_p L_g)(s)=(L_f\# L_g)(s,p)$.
Then $(L_f \#_p L_g)(s)$ extends to an entire function for every $p\in\N$.
Next we fix $p=l-1$, where $l$ is the weight of $g$. Then we can obtain a functional
equation for a related Dirichlet series, which though involves 
``lower-order" $L$-functions
$L_f,$ $L_g$.

Specifically we consider
\begin{eqnarray}
\Lambda_{f, g}(s)=\int_{0}^{\infty} f(ix)\Big (
\int_{x}^{\infty}+\int_{x}^{0}
g(iy)(y-x)^{l-2}dy \Big )x^{s-1}dx
\label{mell1}
\end{eqnarray}

The relation with values of $L_f \#_{l-1} L_g$ is given by
\begin{proposition}
If $L_f(s), L_g(s)$ are the $L$-functions of $f \in S_k, g \in S_l$, then
\begin{multline*}
\Lambda_{f, g}(s)=-2(\frac{1}{2\pi})^{s+l-1}(l-2)!\Gamma(s)(L_f \#_{l-1} L_g)(s)
-\\
(l-2)!(\frac{1}{2\pi})^{s+l-1} \Gamma(s+j)
\sum_{j=0}^{l-2}\frac{(-1)^j}{j!} L_f(s+j)L_g(l-j-1).
\end{multline*}
\end{proposition}
\prf
\begin{multline*}
\Lambda_{f, g}(s)=2 \int_{0}^{\infty} f(ix)
\Big ( \int_{x}^{\infty} g(iy)(y-x)^{l-2}dy \Big )x^{s-1}dx-\\
\int_{0}^{\infty} f(ix)
\Big ( \int_{0}^{\infty} g(iy)(y-x)^{l-2}dy \Big )x^{s-1}dx
\end{multline*}
From  the binomial expansion of $(y-x)^{l-2}$ we deduce
that this equals
\begin{multline*}
-2(\frac{i}{2\pi})^{s+l-1}(l-2)!\Gamma(s)(L_f \# L_g)(s, l-1)  \\
-\sum_{j=0}^{l-2} \binom{l-2}{j}(-1)^j \Big ( \int_0^{\infty}
f(ix)x^{s+j-1} dx \Big ) \Big ( \int_0^{\infty}
g(iy)y^{l-2-j} dy  \Big )
\end{multline*}
This gives the result.
\qed
The analytic continuation of $\Lambda_{f, g}$ can be deduced from this proposition
together with Theorem (\ref{ac}) and the analytic continuation of $\Gamma(s)L_{f}(s)$.
The functional equation is given by
\begin{proposition}
The function $\Lambda_{f, g}(s)$ satisfies
$$\Lambda_{f, g}(s)=w^{(k-l)/2-s+1}i^{s+k+1}\Lambda_{\hat f, \hat g}(k-l-s+2).$$
\end{proposition}
\prf
The change of variables $x \to 1/wx$ implies that $\Lambda_{f, g}(s)$ equals
$$w^{k/2-s} \int_{0}^{\infty} \hat f(ix) \Big (
\int_{1/wx}^{\infty}+\int_{1/wx}^{0} g(iy)(y-\frac{1}{wx})^{l-2}dy)x^{k-s-1}dx
$$
The same change of variables on $y$ gives
$$-i^{-l} w^{(k-l)/2-s+1}\int_{0}^{\infty} \hat f(ix) \Big (
\int_{x}^{\infty}+\int_{x}^{0} \hat g(iy)(x-y)^{l-2}dy)x^{k-s-l+1}dx
$$
which is the desired statement.
\qed

{\small Mathematisches Institut,
Auf der Morgenstelle 10,
72076 T\"ubingen,
Germany,
\tt deitmar@uni-tuebingen.de}

{\small School of Mathematical Sciences,
University of Nottingham,
University Park,
Nottingham
NG7 2RD,
United Kingdom\\
\tt nikolaos.diamantis@maths.nottingham.ac.uk}


\end{document}